# A Survey of directed graphs invariants

Sheng Chen, Yilong Zhang

schen@hit.edu.cn, zhangyl_math@hit.edu.cn

*Dept. of Mathematics, Harbin Institute of Technology, Harbin 150001, P.R.China*

**Abstract.** In this paper, various kinds of invariants of directed graphs are summarized. In the first topic, the invariant $w(G)$ for a directed graph $G$ is introduced, which is primarily defined by S. Chen and X.M. Chen to solve a problem of weak connectedness of tensor product of two directed graphs. Further, we present our recent studies on the invariant $w(G)$ in categorical view.

In the second topic, Homology theory on directed graph is introduced, and we also cast on categorical view of the definition.

The third topic mainly focuses on Laplacians on graphs, including traditional work and latest result of 1-laplacian by K.C.Chang.

Finally, Zeta functions and Graded graphs are introduced, inclduing Bratteli-Vershik diagram, dual graded graphs and differential posets, with some applications in dynamic system.

## Content

# 0. Introduction.

We shall introduce the notation and terminology below on directed graphs. These definitions are primarily prepared for the first topic.

A *directed graph* (abbreviated as *digraph*) $G$ is a pair $G = (V, E)$, with $E \subseteq V \times V$. For an arc $(u, v) \in E$, $u$ is called *initial vertice* of the arc and $v$ is called *end vertice*; *multiple arcs* are arcs with same starting and ending nodes.

The *directed cycle* $C_n$ is the digraph with vertice set $V = \{0, 1, \ldots, n-1\}$ and arc set $\{(0,1), (1,2), \ldots, (n-1, 0)\}$. In particular, $C_1$ consists of a single vertice with a loop. The *directed path* $P_n$ is $C_n$ with the arc $(n-1, 0)$ removed. The *initial vertice* of $P_n$ is $0$ and the *end vertice* of $P_n$ is $n-1$. The *length* of a directed cycle or path is defined to be the number of arcs it contains.

A *route* $R$ of length $n$ in a digraph $G$ is a pair $(S, \sigma)$, where $S = (v_0, v_1, \ldots, v_n)$ is a series of vertices, and $\sigma = (\sigma_1, \ldots, \sigma_n)$ is the orientation of corresponding edges with $\sigma_i = 1$ if $(v_i, v_{i+1}) \in E$, and $\sigma_i = -1$ if $(v_{i+1}, v_i) \in E$, $i = 0, 1, \ldots, n-1$. A route is *closed* if $v_0 = v_n$; a closed route $R$ is *unrepeated* means $v_i \neq v_j$ whenever $i \neq j$, except for $v_0 = v_n$; a route is *open* if $v_0 \neq v_n$. The *weight* of route $R$ is defined $\sigma(R) = \sum_{i=1}^{n}(-1)^{\sigma_i}$. Let $R_1 = (S_1, \sigma_1)$ be a route with $S_1 = (v_0, v_1, \ldots, v_n)$ and $\sigma_1 = (\sigma_1, \ldots, \sigma_n)$; let $R_2 = (S_2, \sigma_2)$ be another route with $S_2 = (u_0, u_1, \ldots, u_m)$ and $\sigma_2 = (\tau_1, \ldots, \tau_n)$. If $v_n = u_0$, the composite route $R_1 R_2 = (S_1 S_2, \sigma_1 \sigma_2)$ is defined to be $S_1 S_2 = (v_0, v_1, \ldots, v_n, u_1, \ldots, u_m)$, $\sigma_1 \sigma_2 = (\sigma_1, \ldots, \sigma_n, \tau_1, \ldots, \tau_n)$.

A *walk* is a route with $\sigma = (1, \ldots, 1)$, there are three kinds of connectedness for digraphs defined as follow. A digraph $G$ is *strong* (strongly connected) if for each pair $u, v \in V$, there is a walk from $u$ to $v$; it is *unilateral* (unilaterally connected) if for each pair $u, v \in V$, there is a walk from $u$ to $v$ or a walk from $v$ to $u$; it is *weak* (weakly connected) if for each $u, v \in V$, there is a route $R(u, v)$ starts at $u$ and ends at $v$. The *diameter* of digraph $d(G)$ is defined to be the maximal weight of routes in $G$, i.e. $d(G) = \max\{w(u, v) | u, v \in V\}$.

Let $G_1, G_2$ be digraphs. Their *tensor product* $G_1 \otimes G_2$ is defined to be the digraph with vertex set $V_1 \times V_2$, and $((u_1, v_1), (u_2, v_2))$ is an arc in $G_1 \otimes G_2$ if $(u_1, v_1)$ and $(u_2, v_2)$ are arcs in $G_1$ and $G_2$.

A digraph map $\psi: G_1 \to G_2$ is a *homomorphism* if $\psi i_1 = i_2 \psi$, and $\psi e_1 = e_2 \psi$, where $i, e: E \to V$ are maps attach each edge to its initial vertice and end vertice, respectively. $\psi$ is said to be *epimorphism* if it is both a homomorphism and a surjevtive.

# 1, The invariant $w(G)$ and weak connectedness of digraphs tensor product

## 1.1 A survey of historic work

Harary [1] raised problems to characterize the digraphs whose tensor product has each of the three kinds of connectedness, i.e. *strong*, *unilateral*, and *weak*. Harary [1] solved unilateral problem in 1966; McAndrew [2] solved strong one in 1963; and Chen [3] solved weak one in 2014.

To solve problem of tensor product with strong connectedness, McAndrew introduced invariant $D(G)$ for digraph $G$, where $D(G) = gcd(n(C))$ is first introduced, where $C$ takes over all of directed cycles of digraph $G = (V, E)$, and $n(C)$ denote the length of cycle $C$. This definition captures the idea of the "distance" between two points in $G$, in the sense of modulo $D(G)$. Further, tensor product with strong connectedness can be characterized by the invariant $D(G)$.

**Theorem 1.1 (McAndrew)** $G_1 \otimes G_2$ *is strong connected if and only if* $G_1$ *and* $G_2$ *are strong, and* $gcd(D(G_1), D(G_2)) = 1$.

To solve problem of weak connectedness, S. Chen and X.M. Chen introduced invariant $w(G)$, where $w(G) = gcd(\sigma(R))$, where $R$ is taking all of the routes of $G$, and the definition of route can be found in Introduction. $w(G)$ is called the *weight* of $G$.

To compare the two definitions, one may notice that $w(G)$ is a generalization of $D(G)$, because a directed cycle is always a route. Therefore $D(G)|w(G)$, if we assume $0|0$ in some cases.

**Theorem 1.2 (S. Chen, X. M. Chen)** *Let both* $G_1$ *and* $G_2$ *are weak connected,*

*(1) If* $w(G_i) > 0$ *(i=1,2),* $G_1 \otimes G_2$ *is weak connected if and only if* $gcd(w(G), w(G')) = 1$ *and one of* $G_i$ *containing sources (sinks) implies other contains no sinks (sources).*

*(2) If one of* $G_i$, *say* $G_1$, *has weight zero, then* $G_1 \otimes G_2$ *is weak if and only if* $G_1$ *is l-chainable, where* $l = d(G_2)$

Remark: $G$ is *l*-chainable if there does *not* exist a pair of permutation matrix $M$, $N$, $A(G)$ being the adjacent matrix of $G$, such that

$$MA(G)N = \begin{bmatrix} A_1 & 0 \\ 0 & A_2 \end{bmatrix}.$$

## 1.2 A categorical view of $w(G)$

Our aim is to characterize the weight $w(G)$ through digraph homomorphism. Our result is: for digraph $G$ with positive weight, there is an epimorphism from $G$ to $C_{w(G)}$; for digraph $G$ with positive weight, there is an epimorphism from $G$ to $P_{d(G)}$. Further, there is a functor from category of digraphs with positive w(G) to category of directed circles, there is a functor from category of digraphs with zero $w(G)$ to category of directed path.

Before coming to weight of digraph, we first study some properties of weight of routes $\sigma(R)$. In this section, we *always assume* that a digraph is weakly connected and has no multiple edges, and the vertice set $V$ is finite.

**Lemma 1.4** *Let $R_1 R_2$ to be the composition of routes, $R^{-1}$ be the inverse of route, and $\varphi$ be the digraph homomorphism, one has*

*(1) $\sigma(R_1 R_2) = \sigma(R_1) + \sigma(R_2)$;*

*(2) $\sigma(R^{-1}) = -\sigma(R)$;*

*(3) $\sigma(\varphi(R)) = \sigma(R)$.*

The formula (1) implies $\sigma$ is additive; and formula (2) implies if the weight of $R$ is negative, its inverse $R^{-1}$ is positive. As $R$ and $R^{-1}$ serve the same role in calculation, we always assume that the weight of $R$ is *positive*. Formula (3) is true since $\varphi(R) = (\varphi(S), \varphi(\sigma))$, where $\varphi(S) = (\varphi(v_0), \varphi(v_1), \ldots, \varphi(v_n))$, and $\varphi(\sigma) = \sigma$.

Next, we reformulate the definition of weight as follow.

**Definition 1.3** *The weight $w(G) = gcd(W)$ is defined for each digraph $G$, where $W = \{n \in Z_+ | \exists closed\ R, \sigma(R) = n\}$, i.e. n belongs to $W$ if there is a closed route R in G such that $\sigma(R) = n$. If $W$ is an empty set, we define $w(G) = 0$.*

An immediate question rises up: How to *calculate* $w(G)$ for a certain digraph? In Definition2, the greatest common divisor is taken over weight of all closed route in $G$, there can be many of routes. Actually, if $W$ is non-empty, it is *infinite*. To see that, if there is a closed route $R$, s.t $\sigma(R) = k > 0$, then $\sigma(R^n) = nk \in W$, for each $n \in Z_+$, therefore $W$ contains an infinite set and $W$ itself is infinite.

However, one notes that $gcd\{\sigma(R_1), \sigma(R_1 R_2)\} = gcd\{\sigma(R_1), \sigma(R_1) + \sigma(R_2)\}$ $= gcd\{\sigma(R_1), \sigma(R_2)\}$, composed routes do *NOT* provide independent information of weights apart from routes $R_1$ and $R_2$; so they could be excluded from the calculation. Now we are ready to present a theorem for calculating $w(G)$.

**Theorem 1.5 (Theorem of Calculating $w(G)$) (Y.L.Zhang)** *$w(G) = gcd(W')$, where $W' = \{n \in Z_+ | \exists close\ and\ unrepeated\ R, \sigma(R) = n\}$, Since the set of vertices $V(G)$ is finite, $W'$ is finite, therefore $w(G)$ can be calculated(through we have not yet found a combinatorial way).*

The following to theorem is our first characterization of $w(G)$.

**Theorem 1.6 (Y.L.Zhang)** *For a digraph G with $w(G) > 0$. There is an epimorphism $\varphi: G \to C_{w(G)}$. Conversely, let S be the set of positive integer n such* *that there is an epimorphism from G to $C_n$, then S is a finite set, and the largest integer of S equal to $w(G)$.*

Proof: (1) We shall denote the vertices of $C_{w(G)}$ as $\{0,1,2,\ldots,w(G)-1\}$. Suppose $w(G) > 0$ for G, then for an arbitrarily fixed vertice $p \in V(G)$, set $\varphi(p) = 0$, and for $q \in V(G)$, define

$$\varphi: V(G) \to V(C_{w(G)})$$
$$q \to \sigma(R(p,q)) \mod(w(G))$$

where $R(p,q)$ is an arbitrary route start at $p$ and end at $q$, this is well-defined since for another route $R'(p,q)$. Since $R^{-1}(q,p)R'(p,q)$ is a closed route in G, so $w(G)|\sigma(R^{-1}(q,p)R'(p,q))$. We have $\sigma(R'(p,q)) = \sigma(R(p,q)) + \sigma(R^{-1}(q,p)) + \sigma(R'(p,q)) = \sigma(R(p,q)) + \sigma(R^{-1}(q,p)R'(p,q)) \equiv \sigma(R(p,q)) \mod(w(G))$.

The behavior of $\varphi$ on edges $E(G)$ is uniquely determined by its behavior on $V(G)$, therefore $\varphi$ is a homomorphism. Since $\varphi$ is obviously a surjective, $\varphi$ is indeed an epimorphism.

(2) We shall only prove that for $n > w(G)$, there cannot be a homomorphism from G to $C_n$. If it is possible, and $\rho: G \to C_n$ is a digraph homomorphism, by definition of weight, there exists a route R in G, such that n is not a divisor of $\sigma(R)$. We shall show that the induced homomorphism $\rho|_R: R \to C_n$ leads to contradiction.

To do this, we shall induce *collapse* of route. Graphically,

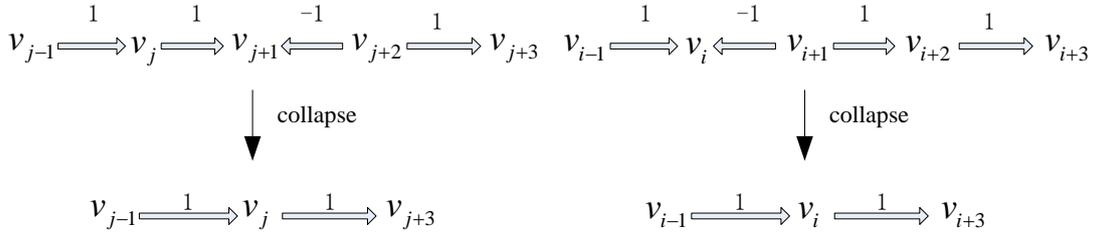

Fig 1.1 the collapse of route

Formally, in a closed route R, if there is a pair of adjacent edges has opposite orientations, i.e. $\exists i$, $\sigma_i = 1$ and $\sigma_{i+1} = -1$, or $\sigma_i = -1$ and $\sigma_{i+1} = 1$, a collapse of R is a closed route $R_1$ by deleting edges $v_i v_{i+1}$ and $v_{i+2} v_{i+1}$ from R, or deleting edges $v_{i+1} v_i$ and $v_{i+1} v_{i+2}$ according to the orientations; deleting $v_{i+1}$; then identify $v_i$ and $v_{i+2}$.

Collapse induces digraph homomorphism $\rho_1: R_1 \to C_n$. We repeat collapse recursively and the process shall terminate at step k with $R_k = C_{\sigma(R)}$, for some integer k, since R is of finite length and its weight is invariant under the collapse, i.e. $\sigma(R) = \sigma(R_1) = \cdots = \sigma(R_k)$.

Therefore, $\rho_k: C_{\sigma(R)}(= R_k) \to C_n$ is a digraph homomorphism. However, n is not a divisor of $\sigma(R)$, this is impossible. □

**Theorem 1.7 (Y.L.Zhang)** *Suppose $w(G) = 0$ for a digraph G. There is a unique epimorphism $\varphi: G \to P_{d(G)}$.*

Proof: let $p$ be the initial point of the path with largest distance in $G$, set $\varphi(p) = 0$ be the initial point of $P_m$, where we denote the vertices of $C_{w(G)}$ as $\{0,1,2,\ldots,P_m - 1\}$. Suppose $w(G) > 0$ for $G$, then for an arbitrarily fixed vertice $p \in V(G)$, set $\varphi(p) = 0$, and for $q \in V(G)$, define

$$\varphi: V(G) \to V(P_m)$$
$$q \to \sigma(R(p,q))$$

where $R(p,q)$ is an arbitrary route start at $p$ and end at $q$, this is well-defined, since for another route $R'(p,q)$. Since $R^{-1}(q,p)R'(p,q)$ is a closed route in G, so $\sigma(R^{-1}(q,p)R'(p,q)) = 0$. We have $\sigma(R'(p,q)) = \sigma(R(p,q)) + \sigma(R^{-1}(q,p)) + \sigma(R'(p,q)) = \sigma(R(p,q)) + \sigma(R^{-1}(q,p)R'(p,q)) \equiv \sigma(R(p,q))$.

In the same way, $\varphi$ is a homomorphism, and a surjective. Therefore, $\varphi$ is indeed an epimorphism. The uniqueness is obvious. □

**Lemma 1.8**

(1) *if $\varphi: C_n \to C_m$ is a homomorphism, then $m|n$, and $\varphi_1$ is an epimorhism;*

(2) *if $\epsilon: P_n \to P_m$ is a homomorphism, then $m \geq n$, and if $m > n$, $\varphi_2$ is a monomorphism;*

(3) *if $\theta: P_n \to C_m$ is a homomorphism, then $n \geq m$, and if $n = m$, $\varphi_3$ maps edges isomorphically and identify the initial and terminal point of $P_n$;*

(4) *there is no homomorphism from $C_n$ to $P_m$.*

**Corollary 1.9**

(1) *If $\varphi: R \to R'$ is a homomorphism between two closed routes, then $\sigma(R')|\sigma(R)$.*

(2) *If $\epsilon: R \to R'$ is a homomorphism between two open routes, then $\sigma(R') > \sigma(R')$.*

Proof: This proof is similar to the Theorem1's, we shall only prove (1) here. Suppose $\sigma(R')$ is not a divisor of $\sigma(R)$. Let $R'_1$ to be the collapsed closed route of $R'$, and $R_1$ to be the induced closed route by collapsing corresponding edges and vertices as well as identifying vertices in the natural way, so that $\varphi|_{R_1}: R_1 \to R'_1$ is a homomorphism. Repeat collapse finite times, one will get $R_k = C_{\sigma(R)}$, and $R'_k = C_{\sigma(R')}$, since weights are invariant under collapse operation. Therefore

$\varphi|_{R_k}: C_{\sigma(R)} \to C_{\sigma(R')}$ is a homomorphism, which is contradictory to Lemma2. □

**Corollary 1.10**

(1) $w(G_i) > 0$ (i=1,2). If $\psi: G_1 \to G_2$ is an epimorphism, then $w(G_2)|w(G_1)$.

(2) $w(G_i) = 0$ (i=1,2). If $\varepsilon: G_1 \to G_2$ is an homomorphism, then $d(G_2) > d(G_1)$.

Proof: (1) If $w(G_2)$ is not a divisor of $w(G_1)$, there exists a closed route $R_0$ in $G_2$, such that for any route $R$ in $G_1$ $\sigma(R_0)$ does not divide $\sigma(R)$, this is contradictory to Corollary1 (1), since $\psi$ is a epimorphism. Proof of part (2) is similar. □

**Main Theorem 1.11 (Y.L.Zhang)**

(1) $w(G_i) > 0$ (i=1,2). If $\psi: G_1 \to G_2$ is an epimorphism, and then $\psi$ induce an epimorphism, $\tilde{\psi}: C_{w(G_1)} \to C_{w(G_2)}$, making the diagram1 commutes. In other words, let *Digraph$_{w>0}$* to be the category of directed graphs with positive weights, *C* to be the category of directed circles, then **Φ**: *Digraph$_{w>0}$* → *C* is a functor.

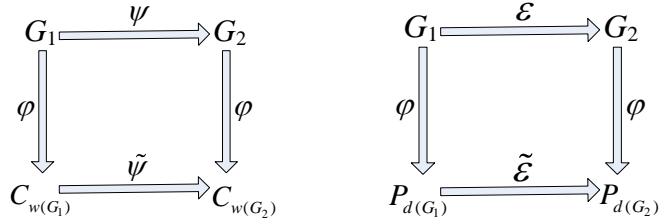

Diagram1,         Diagram2

*Fig1.2 Commutative diagrams for $\varphi$*

(2) $w(G_i) = 0$ (i=1,2). If $\varepsilon: G_1 \to G_2$ is an homomorphism, and then $\varepsilon$ induce an epimorphism $\tilde{\varepsilon}: P_{d(G_1)} \to P_{d(G_2)}$, making the diagram2 commutes. In other words, let *Digraph$_{w=0}$* to be the category of directed graphs with zero weight, *P* to be the category of directed circles, then **Φ**: *Digraph$_{w=0}$* → *P* is a functor.

Proof: (1) We shall define

$$\tilde{\psi}: C_{w(G_1)} \to C_{w(G_2)}$$
$$\varphi(p) \to \varphi \circ \psi(p)$$

First, $\tilde{\psi}$ is well-defined: take $\varphi(p) \in C_{w(G_1)}$, if $q \in G_1$ and $\varphi(q) = \varphi(p)$, then $\exists k \in Z$, for $\forall R(p,q)$ in $G_1$, $\sigma(R(p,q)) = kw(G_1)$. Since $\psi$ is a digraph homomorphism, by Lemma1, one has $\sigma\big(\psi(R(p,q))\big) = \sigma(R(p,q)) = kw(G_1)$, and $\psi(R(p,q))$ is a route in $G_2$ which starts at $\psi(p)$, and ends at $\psi(q)$. Therefore, $w(G_2) | \psi(R(p,q))$. We have $\varphi \circ \psi(p) = \varphi \circ \psi(q)$, so $\tilde{\psi}$ is well-defined.

Secondly, $\tilde{\psi}$ is a homomorphism. To prove this, one note that $\tilde{\psi}(\varphi(p) + 1) = \tilde{\psi}(\varphi(p)) + 1$, therefore $\tilde{\psi}$ maps edges to edges, and is indeed a homomorphism, and then obviously a epimorphism. □

## 2. Homology and Cohomology on digraphs

Grigor'yan, Lin, Muranov, and S-T. Yau [4~6] systematically introduce the *Homology* and *Cohomology* theory on finite digraphs, based on the theories of *p-path* and *p-form* on finite sets, respectively. The boundary operator $\partial$ acting from p-path to (*p*-1)-path and derivative operator $d$ acting from p-form to (*p*+1)-form are naturally defined, satisfying $d^2 = 0, \partial^2 = 0$.

In this section, we briefly introduce Homology theory on digraphs, and present our result based on the definition. As for the Cohomology theory on digraphs, it is dual to Homology case.

### 2.1 Homology on directed graphs

Let *K* be a fixed field, and $G = (V, E)$ be a digraph with the finite set *V* of vertices. An elementary *p*-path on *V* is any (ordered) sequence $i_0, ..., i_p$ of *p*+1 of vertices that will be denoted by $e_{i_0,...,i_p}$. Clearly, $\{e_{i_0,...,i_p}\}$ is a basis. Consider the free *K*-module $\Lambda_p = \Lambda_p(V)$ which is generated by elementary p-paths $e_{i_0,...,i_p}$, whose elements are called *p*-paths, therefore each *p*-path has the form

$$v = \sum_{i_0,...,i_p} v^{i_0,...,i_p} e_{i_0,...,i_p}$$

**Definition 2.1** *The boundary operator* $\partial: \Lambda_p \to \Lambda_{p-1}$ *is defined by*

$$\partial e_{i_0,...,i_p} = \sum_{q=0}^{p} (-1)^q e_{i_0...\hat{i_q}...i_p} \qquad (2.1)$$

*where* $e_{i_0,...,i_p}$ *is elementary p-path.*

Note that the Definition 1. is equivalent to:

$$(\partial v)^{i_0,\ldots,i_{p-1}} = \sum_{k\in V}\sum_{q=0}^{p}(-1)^q v^{i_0\ldots i_{q-1}k i_q\ldots i_{p-1}}.$$

The product of any two paths $u \in \Lambda_p$ and $v \in \Lambda_q$ is defined as $uv \in \Lambda_{p+q+1}$ as follow:

$$(uv)^{i_0\ldots i_p j_0\ldots j_q} = u^{i_0\ldots i_p} v^{j_0\ldots j_q}. \qquad (2.2)$$

For example, if $u = e_{i_0,\ldots,i_p}$ and $v = e_{j_0,\ldots,j_q}$, then

$$e_{i_0,\ldots,i_p} e_{j_0,\ldots,j_q} = e_{i_0,\ldots,i_p j_0,\ldots,j_q}.$$

Based on Definition 1, one can easily observe that

**Theorem 2.2** $\partial^2 = 0$.

**Theorem 2.3** *For $u \in \Lambda_p$ and $v \in \Lambda_q$, we have*

$$\partial(uv) = (\partial u)v + (-1)^{p+1}(\partial v). \qquad (2.3)$$

These basic results come parallelly from the Simplicial homology theories, and the only difference is that the bases are indexed by vertices set.

An elementary p-path $e_{i_0,\ldots,i_p}$ is called *regular* if $i_k \neq i_{k+1}$ for all k, and let $R_p$ denote the submodule generated by regular p-path. Then a regular path is *allowed* if $i_k i_{k+1} \in E, \forall k = 0,\ldots, p-1$. Let $A_p$ denote the submodule generated by all allowed p-paths, one has $A_p \leq R_p \leq \Lambda_p$.

However, the $\partial$ operator is not necessarily invariant on $A_p$, for example in Fig2, $e_{124} \in \Lambda_2$, but $\partial e_{124} = e_{12} - e_{14} + e_{24} \notin A_1$ since 14 is not an edge of the digraph. Therefore, to construct homology on digraph, we need to consider a submodule of $A_p$ which is $\partial$-invariant.

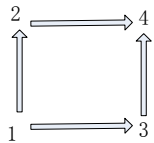

Fig2.1 path $e_{124}$ is not $\partial$-invariant

Consider the submodule $\Omega_p = \Omega_p = \{v \in A_p, \partial v \in A_{p-1}\}$. By definition, $\Omega_p$ is $\partial$-invariant and its elements are called $\partial$-invariant p-paths, and by Theorem 1, one has well-defined chain-complex:

$$\cdots \xrightarrow{\partial_{p+1}} \Omega_p \xrightarrow{\partial_p} \Omega_{p-1} \xrightarrow{\partial_{p-1}} \cdots \xrightarrow{\partial_2} \Omega_1 \xrightarrow{\partial_1} \Omega_0 \xrightarrow{\partial_0} 0.$$

It follows that

**Definition 2.4** *The homology of a digraph G=(V,E) is defined to be :*

$$H_p(G) = Ker\partial_p / Im\partial_{p+1} \tag{2.4}$$

Next, a question rises up: which sort of digraphs have non-trivial homology groups, and what is the necessary condition to have non-trivial homology groups ?

The simplest and most important examples of digraphs are *triangle* and *square*, with orientation defined in Fig3.

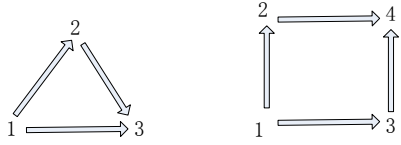

Fig2.2.(1) triangle      (2) square

The first interesting result is that they are the only closed *route* with trivial $H_1$. Here closed route is a digraph with series of distinct vertices $(x_0, x_1, \ldots, x_n)$, with $x_i x_{i+1} \in E$ or $x_{i+1} x_i \in E$. We rephrase the result at hear:

**Theorem 2.5** *Given a close route R, then either $H_1(R) = K$ or $H_1(R) = 0$. And if $H_1(R) = 0$, R is triangle or square.*

To see why $H_1$s of triangle and square are trivial, one notice that $e_{12} - e_{13} + e_{23}$ is closed, that is $\partial(e_{12} - e_{13} + e_{23}) = 0$, and $\partial e_{123} = e_{12} - e_{13} + e_{23}$, therefore $e_{123} \in \Omega_2$, and $\Omega_2 = \text{span}\{e_{123}\}, \Omega_1 = \text{span}\{e_{12} - e_{13} + e_{23}\}, \Omega_2 = \Omega_1 \cong K$, so $H_1(triangle) = 0$. As for square, notice that $\partial(e_{124} - e_{134}) = e_{12} + e_{14} - e_{13} - e_{34}$. In the same way that $\Omega_2 = \Omega_1 \cong K$, so $H_1(square) = 0$.

We end with another theorem which describes in homology theory on digraphs, "interesting" digraph must contain a triangle or a square

**Theorem 2.6** *if a digraph has nontrivial $H_p$ (p>2), it must contains a triangle or a square as sub-digraph.*

## 2.2 Remark: A categorical view.

**Theorem 2.7 (Y.L. Zhang)** *The formula (2.4) defines a functor from the category of digraphs to the category of abelian groups.*

Proof: Let $\psi: G_1 \to G_2$ be digraph homomorphism, given $e^{(1)}_{i_0,\ldots,i_p} \in A^{(1)}_p$, $(i_k, i_{k+1}) \in E^{(1)}$, then $(\varphi(i_k), \varphi(i_{k+1})) \in E^{(2)}$, for $k = 0, \ldots n-1$, therefore $\psi(e^{(1)}_{i_0,\ldots,i_p}) = e^{(2)}_{\psi(i_0),\ldots,\psi(i_p)} \in A^{(2)}_p$. Let $\psi$ still denotes the restriction of $\psi$ on $\Omega_p$. Therefore $\partial(e^{(1)}_{i_0,\ldots,i_p}) = \sum_{q=0}^{p}(-1)^q e^{(1)}_{i_0\ldots\hat{i_q}\ldots i_p}$, and $\psi \circ \partial(e^{(1)}_{i_0,\ldots,i_p}) = \psi(\sum_{q=0}^{p}(-1)^q e^{(1)}_{i_0\ldots\hat{i_q}\ldots i_p}) = \sum_{q=0}^{p}(-1)^q e^{(2)}_{\psi(i_0)\ldots\psi(\hat{i_q})\ldots\psi(i_p)}$. On the other hand, $\psi(e^{(1)}_{i_0,\ldots,i_p}) = e^{(2)}_{\psi(i_0)\ldots\psi(i_p)}$, and $\partial \circ \psi(e^{(1)}_{i_0,\ldots,i_p}) = \partial(e^{(2)}_{\psi(i_0)\ldots\psi(i_p)}) = \sum_{q=0}^{p}(-1)^q e^{(2)}_{\psi(i_0)\ldots\widehat{\psi(i_q)}\ldots\psi(i_p)}$. One notes that $e^{(2)}_{\psi(i_0)\ldots\psi(\hat{i_q})\ldots\psi(i_p)} = e^{(2)}_{\psi(i_0)\ldots\widehat{\psi(i_q)}\ldots\psi(i_p)}$, so we have proved $\psi \circ \partial = \partial \circ \psi$, i.e. the diagram

$$\begin{array}{ccc} \Omega^{(1)}_p & \xrightarrow{\partial} & \Omega^{(1)}_{p-1} \\ \psi \downarrow & & \downarrow \psi \\ \Omega^{(2)}_p & \xrightarrow{\partial} & \Omega^{(2)}_{p-1} \end{array}$$

Fig. 2.3 the commutative diagram of boundary operator $\partial$ and digraph homomorphism $\psi$

commutes, and $\psi: G_1 \to G_2$ induced $\tilde{\psi}: H_p(G_1) \to H_p(G_2)$. To conclude, $\Psi$ is indeed a functor from the category of digraphs to the category of abelian groups. □

## 3. Laplacians on graphs

We shall introduce theories of Laplacians on (undirected) graphs.

The common definition of *Laplacians* on a connected graph $G = (V, E)$ is $L(G) = DD^T$, where $D$ is the incidence matrix of $G$. Therefore $L(G)$ has size $n \times n$, and $|V| = n$. An easy observation shows $L(G) = \Delta(G) - A(G)$, where $\Delta(G)$ is a diagonal matrix with *ii*-entry filled by valency of vertice $i$, $A(G)$ is adjacent matrix of $G$. One can see such definition in [7] and [8], another definition of Laplacian on directed graphs could be found in [17] and [18], we will also introduce as follow.

We first summarize classical results, and then have a brief view of K.C.Chang's work of 1-laplacians on graphs.

## 3.1 A survey of classical results

**I. Basic theorems**

Here we list some basic results related to $L(G)$, one can find further details in [9] and [10].

**Lemma 3.1** rk($L(G)$)=$n$-1, and the ker$L(G)$ is the cyclic group generated by the vector $(1,1,\ldots,1)^T$ in $\mathbf{Z}^n$.

**Lemma 3.2** $L(G)$ is a real symmetric matrix therefore admits n real eigenvalues, with $\lambda_1 < \lambda_2 \leq \cdots \leq \lambda_n$, with $\lambda_1 = 0$ and $\lambda_2$ strict positive.

In the following passage, one will see that much of what theories will concentrate on the information determined by the particular eigenvalue $\lambda_2$.

**Lemma 3.3** Let $x$ be a vector, and we have

$$(x, L(G)x) = x^T L(G) x = \sum_{uv \in E} (x_u - x_v)^2.$$

**II. Potential theory with real coefficient**

Biggs [10] systematically summarized *potential theory on graphs*, we only have a look at the first part here, i.e. potential theory with *real coefficient*.

The theory is based on studying the incidence matrix $D$, as linear mapping between $C^0(G; \mathbf{R})$ and $C^1(G; \mathbf{R})$, the vector spaces of real-valued functions on $V$ and $E$, respectively.

Indeed, interpreting the elements of these spaces as column vectors, we have linear operators

$$D: C^1(G; \mathbf{R}) \to C^0(G; \mathbf{R}), \text{ and } D^T: C^0(G; \mathbf{R}) \to C^1(G; \mathbf{R}).$$

Given a $\alpha \in C^0$, the formula for $D^T \alpha$ is defined as:

$$D^T \alpha(e) = \alpha(h(e)) - \alpha(t(e)).$$

where $h$ means head, and $t$ means tail.

One first note that Given $f \in C^1$, the defnition of $D$ provides an explicit formula for $Df \in C^0$:

$$(Df)(v) = \sum_{h(e)=v} f(e) - \sum_{t(e)=v} f(e).$$

Therefore, a function $f$ is in the kernel of $D$ if and only if for all $v \in V$,

$$\sum_{h(e)=v} f(e) = \sum_{t(e)=v} f(e).$$

This means that the net accumulation of $f$ at every vertice $v$ is zero, and for this reason a function $f \in kerD$ is sometimes known as a *flow*.

We define inner product for $f, g \in C^1(G; \mathbf{R})$: $\langle f, g \rangle = \sum_e f(e)g(e)$. According to the general theory of the vector space, there is a orthogonal decomposition of $C^1$:

$$C^1(G, \mathbf{R}) = \ker D \oplus (\ker D)^\perp.$$

Further, [10] part 1 provide another interpretation of the orthogonal decomposition of $C^1$, as

$$C^1(G, \mathbf{R}) = (flows) \oplus (cuts).$$

where *cut* is interpreted as characteristic functions of boundary of vertice subset $\partial U$, with $U \subseteq V$.

### III. Critical groups

Consider $L(G)$ as an $n \times n$ integer matrix as a map from $\mathbf{Z}^n$ to itself. And there is a natural short-exact sequence:

$$0 \to kerL(G) \xrightarrow{i} \mathbf{Z}^n \xrightarrow{L(G)} \mathbf{Z}^n \xrightarrow{\pi} cokerL(G) \to 0.$$

Since $G$ is connected, The cokernel of $L(G)$ has the form $cokerL(G) = \mathbf{Z}^n/ImL(G) = \mathbf{Z} \oplus K(G)$ where $K(G)$ is defined to be the *critical group* of $G$[11].

It follows from Kirchhoff's Matrix–Tree Theorem that the order $\frac{1}{n}|K(G)| = \kappa(G)$, where $\kappa(G)$ is the number of spanning trees in $G$.

Two integer matrices $A$ and $B$ are equivalent if there exist integer matrices $P$ and $Q$ of determinant $\pm 1$ such that $PAQ = B$. Given a square integer matrix A, its Smith normal form is the unique equivalent diagonal matrix $S(A) = \text{diag}(s_1, s_2, \ldots, s_n)$, whose entries $s_i$ are nonnegative and $s_i$ divides $s_{i+1}$. The $s_i$ and $d_i = \prod_{1 \le j \le i} s_j$ are known as the i*nvariant factors* and *determinantal divisors* of A respectively.

The structure of the critical group is closely related to the Laplacian matrix: if the Smith normal form of $L(G)$ is $\text{diag}(s_1, s_2, \ldots, s_n)$, $K(G)$ is the *torsion* subgroup of $Z_{s_1} \times Z_{s_2} \times \ldots \times Z_{s_n}$.

Research related to Critical groups and Critical ideals can be found in [12~15].

### IV. Cheeger Inequality

The famous Cheeger's inequality from Riemannian Geometry has a discrete analogue involving the Laplacian matrix. This is one of the most important theorems in spectral graph theory. One can find further topics on spectral graph theory in [16].

The *Cheeger constant* of a graph is a numerical measure of whether or not a graph has a "bottleneck". More formally, the Cheeger constant $h(G)$ of a graph $G$ on $n$ vertices is defined as

$$h(G) = \min_{|S| \leq n/2} \frac{|\partial(S)|}{|S|}.$$

where the minimum is over all nonempty sets S of at most n/2 vertices and $\partial(S)$ is the edge boundary of S, i.e., the set of edges with exactly one endpoint in S.

The famous Cheeger Inequality relate the Cheeger constant $h(G)$ and the second eigenvalue $\lambda_2$.

$$\frac{\lambda_2}{2} \leq h(G) \leq (2\lambda_2)^{1/2}. \tag{3.1}$$

## 3.2 Recent work of 1-laplacian on graphs, by K. C. Chang

Traditionally, The Laplace operator is a differential operator acting on functions defined on a manifold $M$, with $\Delta u = div(\nabla u)$. In 2009, K.C.Chang [19] studied nonlinear eigenvalue theory: 1-laplacian operator on $\boldsymbol{R}^n$, where $\Delta_1 u$ is defined as $div(\frac{\nabla u}{|\nabla u|})$, and the eigenvalue problem is to find a pair $(\mu, u) \in \boldsymbol{R}^n \times BV(M)$, satisfying

$$\Delta_1 u \in \mu \operatorname{sgn}(u).$$

Recently(Dec,3,2014), K.C.Chang studied spectrum of 1-laplacian on graphs in [20]. In Chang's way, the 1−Laplace operator on graphs is defined to be:

$$\Delta_1 x = D^T Sgn(Dx).$$

where $D$ is the incidence matrix, and Sgn: $\boldsymbol{R}^n \rightarrow (2^R)^n$ is a set valued mapping, with Sgn(y)=(Sgn($y_1$), Sgn($y_2$),…, Sgn($y_n$)), $\forall$ y = ($y_1, y_2, …, y_n$), in which

$$Sgn(t) = \begin{cases} 1 & \text{if } t > 0 \\ -1 & \text{if } t < 0 \\ [-1,1] & \text{if } t = 0 \end{cases}$$

is a set-valued function.

**Main Theorem. (K.C.Chang)** *Assume that G=(V,E) is connected, then*

$$\mu_2 = h(G). \tag{3.2}$$

Contrasting (3.1) where Cheeger constant and second eigenvalue are related by an inequality for the linear spectral theory, (3.2) provide an identity between them in 1-Laplacian settings. This is the motivation of Chang's study of the nonlinear eigenvalue theory.

# 4. Zeta functions

It is well known that Zeta function plays an important role in number theory and algebraic geometry. It should be noticed that many papers defined and studied various kinds of zeta functions for undirected graphs, directed graphs or symbolic dynamics. For example, for a topological dynamic system $(X,T)$, where $X$ is a topological space, and $T: X \to X$ is a homeomorphism, the Artin-Mazur zeta function $\varsigma_T$ for $(X,T)$ is defined to be

$$\zeta_T(t) = \exp(\sum_{n=1}^{\infty} \frac{p_n}{n} t^n),$$

where

$$p_n = \#\{x \in X : T^n x = x\}.$$

Kim [32] obtained the Artin-Mazur zeta function for a flip system $(X,T,F)$, i.e, $X$ with two homeomorphisms, such that $TF = FT^{-1}$ and $F^2 = id$. This could be interpreted as an infinite dihedral group $D_\infty$ act on $X$. Also, Kim found that when the underlying Z-action is conjugate to a topological Markov shift, the flip system is represented by a pair of matrices, and its zeta function is expressed explicitly in terms of the representation matrices.

For other researches related to this topics, one can refer to [21~31].

# 5. Graded graphs

In this section, we shall introduce graded graphs, differential posets, and related researches, such as Quantized dual graded graphs, Bratteli-Vershik diagram.

## 5.1 Graded graphs and differential posets

A graph is said to be graded if its vertices are divided into levels numbered by integers, so that the endpoints of any edge lie on consecutive levels. Discrete modular lattices and rooted trees are among the typical examples, see Fomin [33].

**Definition 5.1** *A graded graph is a triple $G = (V, \rho, E)$, where*

*(1) V is a discrete set of vertices;*

*(2) $\rho: V \to Z$ is a rank function,*

*(3) E is a multiset of arcs/edges (x,y) where $\rho(y) = \rho(x) + 1$.*

The set $V_n = \{x: \rho(x) = n\}$ is called a level of G. Remark: It is always assumed that each level is finite and $V_0 = \{\hat{0}\}, V_{-1} = V_{-2} = \cdots = \emptyset$.

For a graded graph $G$, one defines an order on $V$ that $a < b$ if and only if $a$ and $b$ lie on consecutive levels and are adjacent, then it could be extended into a partial order. In such way the vertices set becomes a *graded poset* $P = (V, \rho, <)$. Conversely, a graded poset can induce structure of graded graph. We shall interchangeably use the two terms: graded graph and graded poset.

*Differential poset*, introduced by R.P. Stanley [34] is a special kind of graded poset, and is a generalization of Young Lattice.

**Definition 5.2** *A poset P is said to be a differential poset, and in particular to be r-differential (where r is a positive integer), if it satisfies the following conditions:*

*(1) P is graded and locally finite with a unique minimal element $\hat{0}$;*

*(2) for every two distinct elements x, y of P, the number of elements covering both x and y is the same as the number of elements covered by both x and y; and*

*(3) for every element x of P, the number of elements covering x is exactly r more than the number of elements covered by x.*

As an example of differential poset, we have

**Proposition 5.3** *Young's Lattice Y is 1-differential. More generally, $Y^r$ is r-differential.*

Further, in order to study the combinatorial property of differential poset (especially the enumerative property), two linear transform introduced. This work could be found in R.P. Stanley [34]. We shall briefly list some results here.

let $K$ be a field with zero characteristic, and $KP$ be the $K$-vector space with basis $P$. Define two linear transform $U, D: KP \to KP$,

$$Ux = \sum_{y \text{ covers } x} y,$$
$$Dx = \sum_{x \text{ covers } y} y$$

**Theorem 5.4** *Let P be a locally finite graded poset with $\hat{0}$, with finitely many elements of each rank, then P is r-differential poset if and only if $DU - UD = rI$.*

Let $e(x \to y)$ denote the number of shortest non-oriented paths between $x$ and $y$, and $e(y) = e(\hat{0} \to y)$ to be the number of paths going from $\hat{0}$ to $y$. Let

$$\alpha(n \to m) = \sum_{x \in P_n} \sum_{y \in P_m} e(x \to y).$$

In other words, $\alpha(n \to m)$ is the number of paths connecting the nth and mth levels. One can similarly defines $\alpha(n \to m \to p)$. Some enumerative properties are listed as follow:

**Theorem 5.5** *For any graded poset P, one has*

$$\alpha(0 \to n \to 0) = \sum_{x \in P_n} e(x)^2.$$

*Moreover, if P is r-differential poset, one has*

$$\alpha(0 \to n \to 0) = \sum_{x \in P_n} e(x)^2 = r^n n!.$$

For further researches based on the differential posets, R.P. Stanley [35] identify the $P_{[1,2]}$, the poset between level 1 and level 2 of a differential poset $P$ with *hypergraph H*. To do this, one first identify the elements of $P_1$ with the integers $\{1,2,...,r\}$ and element $x_2 \in P_2$ with subset $\{a_1, a_2, ..., a_t\}$ of $V = \{1,2,...,r\}$, where $t \leq r$, and $a_i$ is the element in $P_1$ covered by $x_2$. This construction uniquely corresponds $P_{[1,2]}$ to a hypergraph.

## 5.2 Quatized dual graded graphs, Bratteli-Vershik diagram

Lam [36] studied *Quantized dual graded graphs*, where the operator relation is $DU - qUD = rI$.

*Bratteli-Vershik diagram* was introduced by Ola, Bratteli[37] in 1972 in the theory of operator algebras to describe directed sequences of finite-dimensional algebras: it played an important role in Elliott's classification of AF-algebras and the theory of subfactors. Subsequently Anatoly Vershik associated dynamical systems with infinite paths in such graphs.[38]

As an application, [39] shows that the class of essentially minimal pointed dynamic systems $(X, \varphi, y)$ can be defined and classified in terms of an *ordered Bratteli-Vershik diagram*.

Research Paper 88, 11 pp.

[37] O. Bratteli. Inductive limits of finite dimensional C*-algebras. Trans. Amer. Math. Soc. 171 (1972), 195–234.

[38] A. M. Vershik. A theorem on Markov periodic approximation in ergodic theory, Zap. Nauchn. Sem. Leningrad. Otdel. Mat. Inst. Steklov. (LOMI) 115 (1982), 72–82, 306 .

[39] R.H. Herman; I.F. Putnam; C.F. Skau. Ordered Bratteli diagrams, dimension groups and topological dynamics. Internat. J. Math. 3 (1992), no. 6, 827–864.